\documentclass[11pt, a4paper,leqno]{amsart}
\usepackage{amsmath,amsthm,amscd,amssymb,amsfonts,amsbsy}
\usepackage{latexsym}
\usepackage{exscale}

\usepackage[colorlinks=true, pdfstartview=FitV, linkcolor=blue, citecolor=red, urlcolor=blue, backref=page]{hyperref} 

\usepackage{latexsym}
\usepackage{enumitem}
\usepackage{calc}
\usepackage{tikz}
\usepackage{tkz-euclide}
\usetkzobj{all}
\usetikzlibrary{intersections}
\usetikzlibrary{shapes,arrows,calc}

\calclayout
\allowdisplaybreaks

\theoremstyle{plain}
\newtheorem{theorem}[equation]{Theorem}
\newtheorem{lemma}[equation]{Lemma}

\theoremstyle{definition}
\newtheorem{definition}[equation]{Definition}

\theoremstyle{remark}
\newtheorem{remark}[equation]{Remark}

\numberwithin{equation}{section}
\usepackage{chngcntr} 
\counterwithin{figure}{section}

\renewcommand{\emptyset}{\mbox{\textup{\O}}}

\newcommand\rh{\mathcal{H}}
\newcommand{\hd}{\mathcal{H}-\mbox{dim}\;}
\newcommand{\mc}[1]{\mathcal{#1}}

\newcommand\rd{\mathrm{d}}

\def\XXint#1#2#3{{\setbox0=\hbox{$#1{#2#3}{\int}$ }
\vcenter{\hbox{$#2#3$ }}\kern-.55\wd0}}

 \def\doublespaced{\baselineskip=\normalbaselineskip
    \multiply\baselineskip by 2}
\def\doublespace{\doublespaced}

\doublespace

\newcommand{\de}{\delta}

\newcommand{\la}{\lambda}

\renewcommand*{\backref}[1]{}
\renewcommand*{\backrefalt}[4]{%
 \ifcase #1 (Not cited.)%
   \or        (Cited on page~#2.)%
    \else      (Cited on pages~#2.)%
    \fi}
\begin{document}
\allowdisplaybreaks

\title[On the absolute continuity of p-harmonic measure and $\rh^{n-1}$]{On the absolute continuity of p-harmonic measure and surface measure in Reifenberg flat domains}

\address{Murat Akman \\
               Instituto de Ciencias Matem\' aticas \\ CSIC-UAM-UC3M-UCM, Consejo Superior de Investigaciones Cient{\'\i}ficas\\ Calle de Nicol\' as Cabrera, 13-15, E-28049 Madrid, Spain\\
}
		\email{murat.akman@icmat.es}

\author{Murat Akman}



\begin{abstract}
In this paper, we study the set of absolute continuity of p-harmonic measure, $\mu$, and $(n-1)-$dimensional Hausdorff measure, $\rh^{n-1}$, on locally flat domains in $\mathbb{R}^{n}$, $n\geq 2$. We prove that for fixed $p$ with $2<p<\infty$  there exists a Reifenberg flat domain $\Omega\subset\mathbb{R}^{n}$, $n\geq 2$ with $\rh^{n-1}(\partial\Omega)<\infty$ and a Borel set $K\subset\partial\Omega$ such that $\mu(K)>0=\rh^{n-1}(K)$ where $\mu$ is the p-harmonic measure associated to a positive weak solution to p-Laplace equation  in $\Omega$ with continuous boundary value zero on $\partial\Omega$. We also show that there exists such a domain for which the same result holds when $p$ is fixed with $2-\eta<p<2$ for some $\eta>0$ provided that $n\geq 3$.

This work is a generalization of a recent result of Azzam, Mourgoglou, and Tolsa when the measure $\mu$ is harmonic measure at $x$, $\omega=\omega^{x}$, associated to the Laplace equation, i.e when $p=2$.
\end{abstract}

\maketitle

\tableofcontents

\section{Introduction and statement of main results}
\label{intro}
A well-known result of F. and M. Riesz says that if $\Omega$ is simply connected domain whose boundary has finite length in the plane then harmonic measure and arclength are mutually absolutely continuous. In \cite{M85}, Makarov gives a sharp description of the support of harmonic measure and shows that $\lambda$ which is given below is the proper function to measure the size of the support of $\omega$. In particular, if $\Omega\subset\mathbb{R}^{2}$ is simply connected domain in the plane then $\omega\ll \rh^{\lambda}$ where 
\[
\lambda(r):=r\ \mathrm{exp}\left\{C\sqrt{\log \frac{1}{r} \log\log\log \frac{1}{r}}\right\}
\]
for sufficiently large $C$. Here ``$\ll$'' stands for absolute continuity of the measures and we use``$\perp$'' to denote measures are singular and $\rh^{\lambda}$ denotes the Hausdorff measure with respect to the function $\lambda$ (see \eqref{hausdorffmeasure} for definition of $\rh^{\lambda}$). In \cite{M85}, it is also shown that this result is sharp in the following sense; there is an example of a simply connected domain for which $\omega\perp \rh^{\lambda}$ whenever $C$ is sufficiently small in the definition of $\lambda$. In higher dimensions, due to examples of Ziemer in \cite{Z} and Wu in \cite{Wu}, neither $\rh^{n}|_{\partial\Omega}\ll\omega$ nor $\omega\ll\rh^{n}|_{\partial\Omega}$ are true in general without imposing extra topological or non-topological conditions on $\partial\Omega$. In \cite{DJ}, David and Jerison prove that if $\Omega$ is a non-tangentially accessible (NTA for short and see definition \ref{NTA}) domain and $\partial\Omega$ is Ahlfors-David regular (ADR for short and see definition \ref{ADR}) then harmonic measure is mutually absolutely continuous on $\partial\Omega$ with respect to surface measure, in fact they are $A_{\infty}-$equivalent (see \cite{AHMNT}). In \cite{Ba}, Badger considers the same problem by relaxing ADR property by $\rh^{n-1}(\partial\Omega)<\infty$ and proves that $\rh^{n-1}\ll\omega$ on $\partial\Omega$.  Moreover, he also shows that $\omega\ll\rh^{n-1}\ll \omega$ on the set $A\subset\partial\Omega$ where
\[
A=\left\{x\in\partial\Omega:\, \liminf\limits_{r\to 0} \frac{\rh^{n-1}(\Delta(x,r))}{r^{n-1}}<\infty\right\}.
\] 
Here  $\Delta(x,r)=B(x,r)\cap \partial\Omega$. Badger also conjectures that when $\Omega$ is NTA domain then the same result holds not only on $A\subset\partial\Omega$ but on the whole $\partial\Omega$ (see Conjecture 1.3 in \cite{Ba}). However, it turns out that this is not true in general. In fact, in \cite{AMT}, Azzam, Mourgoglou, and Tolsa construct an example of a Reifenberg flat domain (see definition \ref{RF}) $\Omega$ in $\mathbb{R}^{n}$, $n\geq 3$ with $\rh^{n-1}(\partial\Omega)<\infty$ and a Borel set $E\subset\partial\Omega$ such that 
\[
\omega(E)>0=\rh^{n-1}(E). 
\]

One can consider the same problem for the p-harmonic measure associated with a positive weak solution to p-Laplace equation for $1<p\neq 2<\infty$. To define p-harmonic measure and p-Laplace equation, we let $ \Omega\subset \mathbb{R}^{n}$ be a domain  and let $N$ be a neighborhood of $\partial\Omega$.  Fix  $p$,  $1 < p < \infty$, and suppose that $\hat  u$ is  a positive weak solution  to the p-Laplace equation in $  \Omega\cap N$.   That is,
$ \hat  u \in W^ {1,p} (\Omega\cap N) $ and
\begin{align}   
\label{plaplace}  
\int | \nabla \hat  u |^{p - 2}  \, \langle  \nabla \hat  u , \nabla \theta  \rangle \rd x   = 0    
\end{align}
whenever $ \theta  \in W^{1, p}_0 ( \Omega\cap N)$.  Equivalently, we say that $\hat u $ is p-harmonic in $ \Omega\cap N$.  Observe that if $ \hat  u $ is smooth and $ \nabla \hat  u \not = 0 $
in $\Omega\cap N$ then 
\[
\nabla \cdot ( | \nabla \hat  u |^{ p - 2} \, \nabla \hat  u ) \equiv 0
\]
in the classical sense, where $ \nabla \cdot $ denotes divergence.  We assume that  $ \hat  u $ has  zero  boundary values on $ \partial\Omega \cap  N$ in the  Sobolev sense. More specifically, if $ \zeta \in C^\infty_0 ( \Omega\cap N), $ then
$ \hat  u \, \zeta \in W^{1,p}_0 ( \Omega\cap N). $
Extend $ \hat  u $ to  $N$ by putting
$ \hat  u \equiv 0 $ on $N \setminus \Omega$.   Then  $ \hat  u \in W^{1,p}
(N)$ and it follows from (\ref{plaplace}), as in \cite[Chapter 21]{HKM06}, that
there exists a finite positive Borel measure $ \hat  \mu $  on
$ \mathbb{R}^{n}  $ with support contained in $ \partial\Omega \cap N$  satisfying
\begin{align} 
\label{pharmonic}
 \int | \nabla \hat  u |^{p - 2}  \, \langle  \nabla \hat  u ,
 \nabla \phi  \rangle \rd x = -  \int  \phi  \, \rd \hat  \mu   
\end{align}
 whenever
 $ \phi \in C_0^\infty ( N )$. Existence of $\hat{\mu}$ follows from the maximum principle, basic Caccioppoli inequalities for $\hat{u}$ and the Riesz representation theorem for positive linear functional. We note that if $ \partial\Omega $ is smooth enough and $\nabla u\neq 0$ in $\Omega$, then
\[
\rd \hat  \mu = | \nabla \hat  u |^{p - 1} \, \rd \rh^{n-1}|_{\partial\Omega\cap N}.
\]
 \begin{remark}
 \label{remarkmuisomega}
Note that when $p=2$ in \eqref{plaplace} then we have the usual Laplace's equation. Moreover, if $u$ is the Green's function for Laplace's equation with pole, say $z_{0}\in\Omega$,  then the measure in \eqref{pharmonic} corresponding to this harmonic function $u$ is harmonic measure, $\omega$, relative to $z_{0}$. Note also that, p-Laplace equation in \eqref{plaplace} is degenerate when $p>2$ and is singular when $1<p<2$.  Moreover, nonlinear structure of this pde makes it difficult to work with.
 \end{remark}

We next introduce the notion of the \textit{Hausdorff dimension of a measure}. To this end, let $\hat{r}_{0}>0$ be given, and let $0 <\delta < \hat r_0 $ be fixed. Let $\lambda:[0,\infty)\to[0,\infty)$ be a non-decreasing function with $\lambda(0)=0$. Let $d(\cdot)$ denote the diameter of a set. For a given Borel set $E\subset \mathbb{R}^{n} $, we define $(\de, \la)-$\textit{Hausdorff content} of $E$ in the usual way;
\begin{align*}
 \rh_{ \delta }^\lambda(E):= \inf\left\{\sum\limits_{i} \lambda(d(U_i)):\, \, E\subset \bigcup U_{i}, \, \mbox{each}\, \, U_i\, \, \mbox{is open with}\, \, d(U_{i})<\delta\right\}.
 \end{align*}
Then the \textit{Hausdorff measure} of $E$ is defined by
\begin{align}
\label{hausdorffmeasure}
\rh^\lambda(E):= \lim\limits_{\delta \to 0 } \, \,  \rh_{ \delta }^\lambda(E).
 \end{align}
In case $ \lambda(r) = r^\alpha $ we write $ \rh^\alpha $ for $ \rh^\lambda$.  The Hausdorff dimension of $\hat{\mu}$, denoted by $\hd{\hat{\mu}}$, is defined by
\[
\hd{\hat{\mu}}:=\inf\left\{\alpha:\, \, \exists\, \mathrm{Borel}\, \, E\subset\partial \Omega\, \, \mathrm{s.t.}\, \, \rh^{\alpha}(E)=0\, \, \mathrm{and}\, \, \hat{\mu}(\mathbb{R}^{n}\setminus E)=0\right\}.
\]

We return to our study of singular sets of p-harmonic measure with respect to $\rh^{n-1}$ measure on the boundary of certain domains. The natural candidates, i.e, snowflake type domains which give sharpness  in the harmonic case shown by Makarov, do not provide sharpness as it is observed in \cite{BL05} under the p-harmonic setting. On the other hand, result of David and Jerison described above for harmonic measure is extended to p-harmonic setting for $1<p\neq 2<\infty$ by Lewis and Nystr\"om in \cite{LN11}. To state this result, we let $\Omega\subset\mathbb{R}^{n}$ be a bounded NTA domain with constants $M, r_{0}$ whose boundary is ADR. Let $u$ be p-harmonic in $\Omega\cap B(w,4r)$, $w\in\partial\Omega$, $0<r<r_{0}$ and continuous in $\bar{\Omega}\cap B(w,4r)$ with $u\equiv 0$ on $\Delta(w,4r)$. Extend $u$ to $B(w,4r)$ by defining $u\equiv 0$ on $B(w,4r)\setminus \Omega$ and let $\mu$ be the p-harmonic measure as in (\ref{pharmonic}) associated with $u$. Then it is shown in \cite[Proposition 3.4]{LN11} that $\mu\ll\rh^{n-1}\ll \mu$ on $\partial\Omega$, in fact they are $A_{\infty}-$equivalent. Moreover, it also is proven in \cite{LN11} that Badger's result holds under the p-harmonic setting; if $\Omega$ is NTA domain then $\mu\ll\rh^{n-1}\ll \mu$ on the set $A'\subset \Delta(w,4r)\subset\partial\Omega$ where
\[
A'=\left\{x\in\Delta(w, 4r):\, \liminf\limits_{\rho\to 0} \frac{\rh^{n-1}(\Delta(x,\rho))}{\rho^{n-1}}<\infty\right\}.
\]

The first main result proved in this paper is that there are examples of domains for which absolute continuity of p-harmonic measure and $(n-1)-$dimensional Hausdorff measure does not hold when the domain is NTA, not even locally flat in the sense of Reifenberg. 
\begin{theorem}
\label{main}
Let $p$ be fixed with $2<p<\infty$ and $n\geq 2$. Then there exist a $(\hat{\delta}, \infty)-$Reifenberg flat domain $\Omega\subset\mathbb{R}^{n}$ such that $\rh^{n-1}|_{\partial\Omega}$ is a Radon measure and if $u$ is the p-harmonic function in $\Omega$ with continuous zero boundary value on $\partial\Omega$ and $\mu$ is the p-harmonic measure associated with $u$ as in (\ref{pharmonic}) then there exists a Borel set $K\subset\partial\Omega$ such that  
\[
\mu(K)>0=\rh^{n-1}(K).
\]
\end{theorem}
The second result we obtain in this paper concerns existence of such domains when $p\in (1,2)$. In this case we use a result from \cite{LNV11} to conclude that there exist \textit{Wolff snowflakes} such that the sign of certain integral in the Wolff's program is independent of $p$ when $p$ is in an open interval containing $2$. In order to show such a relation, in \cite{LNV11}, Lewis, Nystr\"{o}m, and Vogel ``perturb" off the $p=2$ case from \cite{W95,LVV05}. Note that results in \cite{W95,LVV05} are valid only when $n\geq 3$. We now state our second result. 
\begin{theorem}
\label{main2}
Let $p$ be fixed, $2-\eta<p<2$ for some $\eta>0$. Then there exists a $(\hat{\delta}, \infty)-$Reifenberg flat domain $\tilde{\Omega}\subset\mathbb{R}^{n}$, $n\geq 3$, such that $\rh^{n-1}|_{\partial\tilde{\Omega}}$ is a Radon measure. Moreover, if $\tilde{u}$ is the p-harmonic function in $\tilde{\Omega}$ with continuous zero boundary value on $\partial\tilde{\Omega}$ and if $\tilde{\mu}$ is the p-harmonic measure associated with $\tilde{u}$ as in \eqref{pharmonic} then there also exists a Borel set $\tilde{K}\subset\partial\tilde{\Omega}$ such that 
\[
\tilde{\mu}(\tilde{K})>0=\rh^{n-1}(\tilde{K}).
\]
\end{theorem}

As plan of this paper, we state definition of  non-tangentially accessible domains, Reifenberg flatness, Ahlfors-David regularity, and we give some lemmas concerning the regularity of p-harmonic function in NTA domains in section \ref{defn}. We give construction of Wolff snowflake in section \ref{wolff}. Following \cite{AMT} we construct ``enlarged domain $\Omega_{\epsilon}^{+}$'' from certain domain $\Omega$ and using some results from \cite{LNV11} concerning the dimension of p-harmonic measure, we give a proof of Theorem \ref{main} in section \ref{proofofmain}. In section \ref{proofofmain2}, we give a short description of Wolff's program from \cite{W95,LVV05} to construct \textit{Wolff snowflakes} with certain properties. Then we make some observations when $p$ fixed is in an open interval containing $2$ and give a proof of Theorem \ref{main2}. 
\section{Definitions and preparatory lemmas}
\label{defn}
Some notations and definitions are in order to proceed. In the sequel, $c$ will denote a positive constant $\geq 1$ (not necessarily the same at each occurrence), which may depend only on $p, n$,  unless otherwise stated. In general, $c(a_{1}, . . . , a_{n})$ denotes a positive constant $\geq 1$ which may depend only on $p,n,a_{1},...,a_{n}$ not necessarily the same at each occurrence. 

Let $x=(x_{1}, \ldots, x_{n})$ denote points in $\mathbb{R}^{n}$ and let $\overline{E}=\mbox{cl}(E)$, $\mbox{int} E$, $\partial E$, and $E^{\mathsf{c}}$ be the closure, interior, boundary, and the complement of the set $E\subset\mathbb{R}^{n}$ respectively. Let $\mbox{diam}(E)$ be the diameter of a set $E$. Let $\langle \cdot, \cdot\rangle$ be the usual inner product in $\mathbb{R}^{n}$.  Let $d(E,F)$ denote the usual distance between the sets $E$ and $F$ and let $d_{\rh}(E,F)$ denote the Hausdorff distance between the sets $E$ and $F$ which is defined by;
\[
d_{\rh}(E,F):=\max\left(\sup\{d(E,y);\; y\in F\}, \sup\{d(x,F);\; x\in E\}\right).
\]
Let $B(x, r)$ be the usual open ball centered at $x$ with radius $r > 0$ in $\mathbb{R}^{n}$ and let $\rd x$ denote the Lebesque $n-$measure in $\mathbb{R}^{n}$. Let $\Delta(w,r)=\partial\Omega\cap B(w,r)$. For a given number $t>$ and a cube $Q$, let $l(Q)$ be the side length of $Q$ and  let $tQ$ denote the cube whose side length is $tl(Q)$ with the same center as $Q$.

We state the notion of non-tangentially accessible domain which is initially introduced by Jerison and Kenig in \cite{JK}. 
\begin{definition}[\bf NTA domain]
\label{NTA}
A domain $\Omega$ is called \textit{non-tangentially accessible}(NTA) domain if there exist $M\geq 2$ and $r_{0}$ such that the following are fulfilled:
\begin{enumerate}[labelindent=0pt,itemindent=1em,leftmargin=!]
\setlength\itemsep{1em}

\item[(i)] {\it Corkscrew condition:} for any $w\in\partial\Omega$, $0<r<r_{0}$, there exists $a_{r}(w)\in\Omega$ satisfying
\[
M^{-1}r<|a_{r}(w)-w|<r\, \, \, \mbox{and}\, \, \, M^{-1}r< d(a_{r}(w),\partial\Omega).
\] 
\item[(ii)] $\mathbb{R}^{n}\setminus\overline{\Omega}$ satisfies corkscrew condition.
\item[(iii)] {\it Uniform condition:} if  $w\in\partial\Omega$, $0<r<r_{0}$, and $w_{1},w_{2}\in B(w,r)\cap\Omega$ then there exists a rectifiable curve $\gamma:[0,1]\to \Omega$ with $\gamma(0)=w_{1}$ and $\gamma(1)=w_{2}$ such that\\
(a) $\rh^{1}(\gamma)\leq M|w_{1}-w_{2}|$,\\
(b) $\min\{\rh^{1}(\gamma([0,t])), \rh^{1}([t,1]))\}\leq Md(\gamma(t), \partial\Omega)$.
\end{enumerate}
\end{definition}
\begin{remark}
We use the definition of this notion given in \cite{LN11}. Note that (iii) of definition \ref{NTA} is different but equivalent to the Harnack chain condition given in \cite{JK}. 
\end{remark}
Next we give the definition of Reifenberg flatness from \cite{AMT}.
\begin{definition}[\bf $(\delta,r_{0})$-Reifenberg flat domain]
\label{RF}
Let $\Omega$ be a domain and $r_{0},\delta>0$ with $0<\delta<1/2$. Then $\Omega$ is said to be $(\delta,r_{0})$-\textit{Reifenberg flat} provided that the following two conditions hold.
\begin{enumerate}[labelindent=0pt,itemindent=1em,leftmargin=!]
\setlength\itemsep{1em}
\item[(i)] For every $w\in\partial\Omega$ and every $0<r<r_{0}$ there exists a  a hyperplane $\mc{P}(w,r)$ containing $w$ such that 
\[
d_{\rh}(\Delta(w,r), \mc{P}(w,r)\cap B(w,r))\leq \delta r.
\]
\item[(ii)] For every $x\in\partial\Omega$, one of the connected components of 
\[
B(x,r_{0})\cap \{x\in\mathbb{R}^{n};\, \rd(x, \mc{P}(x,r_{0}))\geq 2\delta r_{0}\}
\]
is contained in $\Omega$ and the other is contained in $\mathbb{R}^{n}\setminus\Omega$.
\end{enumerate} 
\end{definition}
We say that $\Omega$ is $(\delta,\infty)$-Reifenberg flat if it is $(\delta, r_{0})$-Reifenberg flat for every $r_{0}>0$. 
\begin{remark}
An equivalent definition of Reifenberg flatness is given in \cite{LN11} and it is remarked that these two definitions are equivalent(see observation after the Definition 1.2 in \cite{LN11}).
\end{remark}
\begin{definition}[\bf Ahlfors-David regular set]
\label{ADR}
 We say that  $\partial\Omega$ is $n$-dimensional \textit{Ahlfors-David regular}(ADR) if there is some uniform constant $C$ such that
\begin{align*}
\frac{1}{C} r^n \leq \rh^{n}(\Delta(x,r)) \leq C\, r^n,\,\, \forall r\in(0, \mbox{diam}(\Omega)),x \in \partial\Omega,
\end{align*}
\end{definition}
We next give some estimates from \cite{LNV11} when $n\geq 3$ and from \cite{BL05} when $n=2$ given under the p-harmonic settings (see Lemmas 3.2-3.6 in \cite{LNV11} and Lemmas 2.6, 2.7, 2.13, 2.14 in \cite{BL05}). For Lemmas \ref{lemma32}-\ref{lemma36}, let $p$ be fixed with $1<p\neq 2<\infty$.
\begin{lemma}
\label{lemma32}
Let $u$ be a positive p-harmonic function in $B(w,2r)\subset\mathbb{R}^{n}$, $n\geq 3$. Then
\[
r^{p-n} \int\limits_{B(w,r/2)} |\nabla u|^{p}\rd x \leq c\left(\max\limits_{B(w,r)} u\right)^{p}
\]
 and
 \[
 \max\limits_{B(w,r)} u\leq c \min\limits_{B(w,r)} u.
 \]
 Moreover, there exists $\beta=\beta(p,n)\in (0,1)$ such that if $x,y\in B(w,r)$ then
 \[
 |u(x)-u(y)|\leq c\left( \frac{|x-y|}{r}\right)^{\beta}\, \max\limits_{B(w,2r)} u.
 \] 
\end{lemma}
For lemmas \ref{lemma35}-\ref{lemma36} let $\Omega$ be an NTA domain in $\mathbb{R}^{n}$ and let $w\in\partial\Omega$, $0<r<r_{0}$.
\begin{lemma}
\label{lemma35}
Suppose that $u$ is non-negative continuous p-harmonic function in $\bar{\Omega}\cap B(w,4r)$ and $u=0$ on $\Delta(w,4r)$. Extend $u$ to $B(w,4r)$ by defining $u\equiv 0$ on $B(w,4r)\setminus\Omega$. Then $u$ has a representative in $W^{1,p}(B(w,4r))$ with H\"older continuous partial derivatives in $\Omega\cap B(w,4r)$. In particular, there exists $\sigma=\sigma(p,n)\in(0,1]$ such that if $x,y\in B(\hat{w},\hat{r}/2)$ where $B(\hat{w},4\hat{r})\subset\Omega\cap B(w,4r)$ then
\[
\frac{1}{c}|\nabla u(x)-\nabla u(y)|\leq \left(\frac{|x-y|}{\hat{r}}\right)^{\sigma}\, \max\limits_{B(\hat{w},\hat{r})}|\nabla u|\leq \frac{c}{\hat{r}}\left(\frac{|x-y|}{\hat{r}}\right)^{\sigma}\, \max\limits_{B(\hat{w},2\hat{r})} u.
\]
If $\nabla u(\hat{w})\neq 0$ then $u$ is real analytic in a neighborhood of $\hat{w}$.
\end{lemma}
Next lemma gives a relation between the p-harmonic function and the p-harmonic measure.
\begin{lemma}
\label{lemma36}
Suppose that $u$ is non-negative continuous p-harmonic function in $\bar{\Omega}\cap B(w,2r)$ and $u=0$ on $\Delta(w,2r)$.
Extend $u$ to $B(w,2r)$ by defining $u\equiv 0$ on $B(w,2r)\setminus\Omega$. As in (\ref{pharmonic}) Then there exists a unique locally finite positive Borel measure $\mu$ on $\mathbb{R}^{n}$ with support in $\Delta(w,2r)$ such that
\[
\int |\nabla u|^{p-2}\langle \nabla u,\nabla \theta\rangle\rd x=-\int\theta\rd \mu
\]
whenever $\theta\in C^{\infty}_{0}(B(w,2r))$. Moreover, there exists $c=c(p,n,M)\in[1,\infty)$ such that if $\tilde{r}=r/c$ then
\[
c^{-1}r^{p-n}\mu(\Delta(w,\tilde{r}))\leq (u(a_{\tilde{r}}(w)))^{p-1}\leq c\, r^{p-n}\mu(\Delta(w,\tilde{r}/2))
\]
where $a_{\tilde{r}}(w)$ is as in definition \ref{NTA}.
\end{lemma}
\section{Construction of \textit{Wolff Snowflakes}}
\label{wolff}
In this section, following \cite{LNV11} when $n\geq 3$ and \cite{BL05} when $n=2$, we describe the construction of {\it Wolff snowflakes} in $\mathbb{R}^{n}$ which is originally introduced by Wolff in \cite{W95}. To this end, let
\[
\Omega_{0}=\{(x',x_{n}),\, x'\in\mathbb{R}^{n-1}, \, x_{n}>0\}\subset\mathbb{R}^{n}.
\]
Set 
\[
Q(r)=\{x'\in\mathbb{R}^{n-1}; -r/2\leq |x_{i}|\leq r/2, \, \mbox{for}\, 1\leq i\leq n-1\}.
\] 
Then $Q(r)$ is a $(n-1)-$dimensional cube with side length $r$ and centered at $0$.  Let $\phi:\mathbb{R}^{n-1}\to \mathbb{R}$ be a piecewise linear function with support contained in $\{x':\, |x'|<1/2\}$ satisfying 
\begin{align}
\label{nablaphi}
\|\nabla\phi\|_{\infty}\leq \theta_{0}.
\end{align}
For fixed large $N$, define $\psi(x')=N^{-1}\phi(Nx')$. Let $b>0$ be a small constant and let $Q$ be an $(n-1)$ dimensional cube with center $a_{Q}$ and length $l(Q)$ contained in some hyperplane. Let $\mbox{cch}(E)$ denote the closed convex hull. Let $e$ be a unit normal to $Q$ and define
\begin{align*}
P_{Q}=\mbox{cch}(Q\cup \{a_{Q}+bl(Q)e\})\, \, \mbox{and}\, \, \tilde{P}_{Q}=\mbox{int}\, \mbox{cch}(Q\cup \{a_{Q}-bl(Q)e\}).
\end{align*}
We set $e=-e_{n}$ for $Q(1)$.  We also define
\[
\Lambda:=\{x\in P_{Q(1)}\cup \tilde{P}_{Q(1)},\, x_{n}\geq \psi(x)\}\, \, \mbox{and}\, \,  \partial:=\{x\in\mathbb{R}^{n}, \, x'\in Q(1),\, x_{n}=\psi(x')\}.
\]
We assume that $N=N(b,M)$ is so large that 
\[
d(\partial\setminus\partial\Omega_{0}, \partial[P_Q(1)\cup \tilde{P}_{Q(1)}])\geq b/100.
\] 
From the construction, it can be easily seen that $\partial\subset Q(1)\times[-1/2, 1/2]$ consists of a finite number of $(n-1)$ dimensional faces. We fix a Whitney decomposition of each face; we divide each face of $\partial$ into $(n-1)-$dimensional cube $Q$, with side lengths $8^{-k}$, $k=1,2,\ldots,$ and $8^{-k}\approx$ to their distance from the edges of the face they lie on. We next choose a distinguished $(n-2)-$dimensional ``side'' for each $(n-1)-$dimensional cube. 

Suppose $\Omega$ is a domain and $Q\subset\partial\Omega$ is an $(n-1)-$dimensional cube with distinguished side $\gamma$. Let $e$ be a unit normal to $\partial\Omega$ on $Q$ and assume that $P_{Q}\cap \Omega=\emptyset$ and $\tilde{P}\subset \Omega$. We form a new domain $\tilde{\Omega}$ as follows. Let $\mc{T}$ be the conformal affine map, i.e., composition of a translation, dilation, and rotation with $\mc{T}(Q(1))=Q$ which fixes dilation, $\mc{T}(0)=a_{Q}$ which fixes translation, $\mc{T}(\{x\in\partial Q(1):\, x_{1}=1/2\})\, \mbox{and}\, \mc{T}(-e_{n})$ in the direction of $e$ which fixes rotation.
Let  $\Lambda_{Q}=\mc{T}(\Lambda)\, \, \mbox{and}\, \, \partial_{Q}=\mc{T}(\partial)$. Then we define $\tilde{\Omega}$ through the relations 
\[
\tilde{\Omega}\cap(P_{Q}\cup \tilde{P}_{Q})
\]
and
\[
\tilde{\Omega}\setminus (P_{Q}\cup \tilde{P}_{Q})=\Omega\setminus (P_{Q}\cup \tilde{P}_{Q}).
\]
Note that $\partial_{Q}$ inherits from $\partial$ a natural subdivision into Whitney cubes with distinguished sides. This process is called ``adding a blip to $\Omega$ along $Q$''. 

To use the process of ``adding a blip'' to construct a Wolff snowflake $\Omega_{\infty}$, starting from $\Omega_{0}$, we first add blip to $\Omega_{0}$ along $Q(1)$ obtaining a new domain $\Omega_{1}$. We then inherit a subdivision of $\partial\Omega_{1}\cap (P_{Q(1)}\cap \tilde{P}_{Q(1)})$ into Whitney cubes with distinguished sides, together with a finite set of edges $E_{1}$ (the edges of the faces of the graph are not in the Whitney cubes). Let $G_{1}$ be the set of all Whitney cubes in the subdivision. Then $\Omega_{2}$ is obtained from $\Omega_{1}$ by adding blip along each $Q\in G_{1}$. From this process, we inherit a family of cubes $G_{2}\subset\partial\Omega_{2}$ (each with a distinguished side) and a set of edges $E_{2}\subset \partial\Omega_{2}$ of $\sigma-$finite $\rh^{n-2}$ measure. Continuing by induction we get $(\Omega_{m})_{m=n-1}^{\infty}$, $(G_{m})_{m=n-1}^{\infty}$, and $(E_{m})_{m=n-1}^{\infty}$ where 
\[
\partial\Omega_{m}\cap (P_{Q(1)}\cap \tilde{P}_{Q(1)})=E_{m}\cup\bigcup\limits_{Q\in G_{m}} Q\, \, \mbox{for}\, m\geq n-1. 
\]
If $N=N(b,M)$ is large enough, then $\Omega_{m}\to \Omega_{\infty}$ in the Hausdorff distance sense. We call $\Omega_{\infty}$ a {\it Wolff snowflake}. We state a result from \cite{LNV11} which says that \textit{Wolff snowkflakes} are locally flat in the sense of Reifenberg. 
\begin{lemma}[{\cite[Lemma 7.1]{LNV11}}]
\label{wolffisRF}
If $\theta_{0}, N^{-1}$ are small enough, depending only $n$ then the \textit{Wolff snowflake} domain $\Omega_{\infty}$ is $(c\theta_{0},\infty)$-Reifenberg flat where $c=c(n)$. 
\end{lemma}
\section{Proof of Theorem \ref{main}}
\label{proofofmain}
In this section we give a proof of Theorem \ref{main} using some results from \cite{LNV11,AMT}. To this end, let $\Omega_{\infty}$ be a \textit{Wolff snowflake} with constants $\theta_{0}, N$ as described in section \ref{wolff}. For fixed $p$, $1<p\neq 2<\infty$, let $u_{\infty}$ be the unique positive p-harmonic function in $\Omega_{\infty}$ with continuous boundary value zero on $\partial\Omega_{\infty}$ and $|x_{n}-u_{\infty}(x)|\to 0$ uniformly as $|x|\to \infty$. Let $\mu_{\infty}$ be p-harmonic measure associated with $u_{\infty}$ as in (\ref{pharmonic}).  A proof of existence and uniqueness of $u_{\infty}$ can be found in \cite[Lemma 6.1]{LNV11}. Let $\Omega'_{\infty}$ be the restriction of $\Omega_{\infty}$ to $Q(1)\times[-1,1]$ and let $\mu'_{\infty}$ be the restriction of $\mu_{\infty}$ to $(Q(1)\times [-1,1])\cap\partial\Omega_{\infty}$.  In \cite{LNV11}, it is shown that Wolff's program in \cite{W95} can also be made to work under the p-harmonic setting and is observed that certain integral has sign(see integral in \eqref{Iphip} and \cite[section 6]{LNV11} for more details). The following lemma can be easily deduced by combining Lemma 7.4 and Proposition 7.6 from \cite{LNV11} when $n\geq 3$ and combining Lemma 3.23 and Theorem 1 from \cite{BL05} when $n=2$.
\begin{lemma}
\label{dlessn-1}
Let $p$ be fixed, $2<p<\infty$, and let $\Omega'_{\infty}$ and $\mu'_{\infty}$ be described as above. Then for some $d>0$ we have
\[
\lim\limits_{r\to 0} \frac{\log \mu'_{\infty}(\Delta(x,r))}{\log r}\leq d<n-1\, \, \mbox{for all}\, \, x\in\partial\Omega'_{\infty}\setminus\Lambda
\]
where $\Lambda\subset\partial\Omega'_{\infty}$ with $\mu'_{\infty}(\Lambda)=0$. Moreover, $\hd{\mu'_{\infty}}\leq d<n-1$.
\end{lemma}
We are now ready to prove Theorem \ref{main}. Under the p-harmonic setting, we closely follow the arguments given in \cite{AMT} after Theorem 4.3. We first observe from Lemma \ref{dlessn-1}, more specifically from the fact $\hd{\mu'_{\infty}}\leq d<n-1$, and the definition of Hausdorff dimension of p-harmonic measure that there is a Borel set $E\subset\partial\Omega'_{\infty}$ such that $\mu'_{\infty}(\mathbb{R}^{n}\setminus E)=0$ and $\rh^{d}(E)=0$. From this observation and once again from lemma \ref{dlessn-1} we also have
\begin{align}
\label{E}
\lim\limits_{r\to 0} \frac{\log \mu'_{\infty}(B(x,r))}{\log r}\leq d<n-1,\, \, \forall x\in E. 
\end{align}
Note that $\Omega'_{\infty}$ is the restriction of $\Omega_{\infty}$ to $Q(1)\times[-1,1]$, therefore, 
\[
\partial\Omega_{\infty}\setminus \{(x',x_{n})\in\mathbb{R}^{n};\, x_{n}=0\}\subset \partial\Omega'_{\infty}.
\]
For ease of notation we let 
\[
\mathfrak{R}^{n-1}:=\{(x',x_{n})\in\mathbb{R}^{n};\, x'\in\mathbb{R}^{n-1}\, \mbox{and}\, x_{n}=0\}.
\]
From (\ref{E}) it follows that for $\alpha$, $0<\alpha<n-1-d$ one can find small enough $\rho$ such that $\mu'_{\infty}(E_{1})>0$ where
\[
E_{1}=\left\{x\in (E\cap\partial\Omega_{\infty})\setminus \mathfrak{R}^{n-1};\, \frac{\log \mu'_{\infty}(B(x,r))}{\log r}<n-1-\alpha,\, \, \forall r\in (0,\rho] \right\}.
\]
We next fix a point $\zeta_{0}\in E_{1}$. By the regularity of p-harmonic measure we can find $\rho_{0}\in (0,\rho]$ and a compact set $K\subset E_{1}\cap B(\zeta_{0},\rho_{0})$  such that for all $x\in K\, \, \mbox{and}\, r\in(0,\rho_{0})$ with the following property
\[
\mu'_{\infty}(K)>0\, \, \mbox{and}\, \, \mu'_{\infty}(B(x,r))>r^{n-1-\alpha}.
\]
The construction yields that $K\subset \partial\Omega'_{\infty}\cap \partial\Omega_{\infty}$ and $\mbox{cl}(\Omega'_{\infty})\subset\mbox{cl}(\Omega_\infty)$. Then using the fact that the support of  $\mu'_{\infty}$ is contained in $(Q(1)\times [-1,1])\cap\partial\Omega_{\infty}$ we have
\begin{align}
\label{K}
\mu_{\infty}(K)>0\, \, \mbox{and}\, \, \mu_{\infty}(B(x,r)\cap \partial\Omega_{\infty})>r^{n-1-\alpha}
\end{align}
for all $x\in K$ and $r\in(0,\rho_{0})$. 

For a given number $t$, $4\leq t$, and given open set $O\subset\mathbb{R}^{n-1}$ we use $\mc{W}_{t}(O)$ to denote the set of maximal dyadic cubes $Q\subset O$ satisfying $tQ\cap Q^{\mathsf{c}}=\emptyset$. Let $0<\epsilon<1/100$  and let $\mc{I}$ be the family of cubes $Q\in \mc{W}_{\epsilon^{-2}}(K^{\mathsf{c}})$ such that 
\[
Q\cap (Q(1)\times [-1,1])\cap\partial\Omega_{\infty}\neq \emptyset.
\]
Note that
\begin{align*}
l(Q)\approx \epsilon^{2} \mbox{ dist}(Q,K)\, \mbox{for all}\, Q\in \mc{I}\, \, \mbox{and}\, \, \partial\Omega'_{\infty}\setminus K\subset \bigcup\limits_{Q\in \mc{I}} Q. 
\end{align*}
For each $Q\in\mc{I}$, fix some point $z_{Q}\in Q\cap \partial\Omega'_{\infty}$. We then define a new domain $\Omega_{\epsilon}^{+}$ by
\[
\Omega_{\epsilon}^{+}:=\Omega'_{\infty}\cup \left(\bigcup\limits_{Q\in\mc{I}}B_{Q}\right)
\, \, \mbox{where}\, \, B_{Q}=B(z_{Q}, \epsilon\mbox{ dist}(z_{Q},K)).
\]
It is observed in \cite[Lemma 2.2]{AMT} that if $\theta,\epsilon$ in the construction of \textit{Wolff snowflake} in section \ref{wolff} are small enough then $\Omega_{\epsilon}^{+}$ is $(c\epsilon^{1/2}, r_{0})$-Reifenberg flat and $K\subset\partial\Omega_{\epsilon}^{+}$ provided that the original domain $\Omega_{\infty}$ is $(\delta, r_{0})-$Reifenberg flat. Note that from Lemma \ref{wolffisRF} we have that \textit{Wolff snowflake} domain $\Omega_{\infty}$ is $(c\theta_{0},r_{0})-$Reifenberg flat where $r_{0}=\infty$. Therefore if we choose $\theta$ and $\epsilon$ small enough and use Lemma 2.2 from \cite{AMT} we have $\Omega_{\epsilon}^{+}$ is a $(c\epsilon^{1/2}, \infty)$-Reifenberg flat domain satisfying
\begin{align}
\label{omegasubsetepsilon}
K\subset\partial\Omega'_{\infty}\cap \partial\Omega_{\epsilon}^{+}\, \, \mbox{and}\, \, \mbox{cl}(\Omega_{\infty})\subset \mbox{cl}(\Omega_{\epsilon}^{+}).
\end{align}
Let $u_{\epsilon}^{+}$ be a positive p-harmonic function in $\Omega_{\epsilon}^{+}$ with continuous boundary value zero on $\partial\Omega_{\epsilon}^{+}$. Let $\mu_{\epsilon}^{+}$ be the p-harmonic measure associated with $u_{\epsilon}^{+}$ as in (\ref{pharmonic}). From the construction of $\Omega_{\epsilon}^{+}$ we have $u_{\epsilon}^{+}\geq u'_{\infty}$ on $\partial\Omega'_{\infty}$. Then it follows from maximum principle for positive p-harmonic functions and \eqref{omegasubsetepsilon} that $u_{\epsilon}^{+}\geq u'_{\infty}$ in $\Omega'_{\infty}$. This observation, Lemmas \ref{lemma32}-\ref{lemma36} and, (\ref{K}) yield
\begin{align}
\label{muK}
\mu_{\epsilon}^{+}(K)>0\, \, \mbox{and}\, \, \mu_{\epsilon}^{+}(B(x,r))>r^{n-1-\alpha},\, \,\forall  x\in K\, \, \mbox{and}\, \, r\in(0,\rho_{0}).
\end{align}
As $\mu_{\epsilon}^{+}$ is a Radon measure which follows from Lemma \ref{lemma36} and satisfies \eqref{muK} and $\Omega_{\epsilon}^{+}$ is $(\hat{\delta}, \hat{r}_{0})$-Reifenberg flat domain, it follows from \cite[Lemma 3.1]{AMT} that $\rh^{n-1}|_{\partial\Omega_{\epsilon}^{+}}$ is locally finite. Let $\Omega:=\Omega_{\epsilon}^{+}$ be the $(\hat{\delta}, \hat{r}_{0})$-Reifenberg flat domain with locally finite surface measure and let $\mu:=\mu_{\epsilon}^{+}$ be the p-harmonic measure as above. The fact $K\subset\partial\Omega$ is a compact set with $\mu(K)>0=\rh^{n-1}(K)$ and (\ref{muK}) imply that proof of Theorem \ref{main} is now complete.
\qed
\section{Proof of Theorem \ref{main2}}
\label{proofofmain2}
To discuss the matter when $1<p<2$, we give a short description of Wolff's program from \cite{LNV11} and then make some observations. To this end, let $p$ be fixed and define
\[
\hat{\Omega}=\hat{\Omega}(\epsilon):=\{x=(x', x_{n}):\, \, x'\in\mathbb{R}^{n-1}\,\, \mbox{and}\, \, x_{n}>\epsilon\hat{\theta}\}
\]
for some function $\hat{\theta}$ and constant $\epsilon>0$. Let $\hat{u}(\cdot, \epsilon)$ be the positive p-harmonic function in $\hat{\Omega}$ with continuous zero boundary value on $\partial\hat{\Omega}$ and $|x_{n}-\hat{u}(x,\epsilon)|\to 0$ uniformly as $|x|\to\infty$, $x\in\hat{\Omega}$. Let
\begin{align}
\label{Iepsilon}
I=I(\epsilon):=\int\limits_{\partial\hat{\Omega}}|\nabla \hat{u}(x,\epsilon)|^{p-1}\log|\nabla \hat{u}(x,\epsilon)|\, \rd\rh^{n-1}.
\end{align}
When $\hat{\theta}\in C_{0}^{\infty}(\mathbb{R}^{n-1})$ then in \cite{LNV11}, it is shown that
\begin{align}
\label{012der}
I(0)=0,\, \, I'(0)=0,\,\, \mbox{and}\, \, I''(0)=\frac{p-2}{p-1}\int\limits_{\mathbb{R}^{n-1}} |\nabla_{x'}\hat{\theta}|^{2}\rd\rh^{n-1}.
\end{align}
where $\nabla_{x'}$ denotes gradient in $x'\in\mathbb{R}^{n-1}$. One has to show that $I$ has a sign when $\epsilon$ is small in order to give a lower and an upper bound for the Hausdorff dimension of the p-harmonic measure. In fact, the choice of $\hat{\theta}$ determines the sign of $I$ when $\epsilon>0$ is sufficiently small, say $0<\epsilon\leq \epsilon_{0}$ where $\epsilon_{0}=\epsilon_{0}(n,p,\hat{\theta})$. When $p=2$ then $I''(0)=0$ in \eqref{012der} and in order to show that $I$ has a sign one has to calculate $I'''(0)$ and $I^{(4)}(0)$ which are done in \cite{W95,LVV05}. Surprisingly, it is observed in \cite{LNV11} that when $p\neq 2$ then the sign of $I$ depends only on $p$. To simplify calculation, one can approximate $\epsilon\hat{\theta}$ by a piecewise linear function $\phi=\phi(\cdot, \epsilon)$ so that the sign of $I$ is preserved. That is, if we let $\tilde{\Omega}:=\{(x'x_{n}): x'\in\mathbb{R}^{n-1}\,\, \mbox{and}\, \, x_{n}>\phi(x')\}$ and if $\tilde{u}$ is the p-harmonic function in $\tilde{\Omega}$ with continuous boundary value zero on $\partial\tilde{\Omega}$
 and $|x_{n}-\tilde{u}(x)|\to 0$ uniformly as $|x|\to \infty$ for $x\in\tilde{\Omega}$. Then it is shown in \cite{LNV11} that both $I$ and $\tilde{I}$ have same sign where
 \begin{align}
\label{Iphip}
\tilde{I}=\tilde{I}(\phi,p):=\int\limits_{\partial\tilde{\Omega}}|\nabla \tilde{u}|^{p-1}\log|\nabla \tilde{u}|\, \rd\rh^{n-1}.
\end{align} 
We then construct a \textit{Wolff snowflake} as we described in section \ref{wolff} relative to $\psi(x')=N^{-1}\phi(x'N)$ for $x'\in\mathbb{R}^{n-1}$ and we let $\Omega_{1}:=\{(x'x_{n}): x'\in\mathbb{R}^{n-1}\,\, \mbox{and}\, \, x_{n}>\psi(x')\}$. We then repeat the process to obtain $\Omega_{\infty}$. Let $u_{\infty},\mu_{\infty}, \mu'_{\infty}$ be the p-harmonic function, p-harmonic measure, and restriction of the p-harmonic measure as in section \ref{wolff} relative to $\Omega_{\infty}$. 

As observed in \cite{LN11}, Wolff program in \cite{W95} (equivalently) says that
\begin{align}
\label{IpIn}
\left\{
\begin{array}{ll}
\hd{\mu'_{\infty}}>n-1 & \mbox{if}\, \, I<0,\\
\hd{\mu'_{\infty}}<n-1 & \mbox{if}\, \, I>0.
\end{array}
\right.
\end{align}
Note that \eqref{Iphip} and \eqref{IpIn} holds for any $p\in (1,\infty)$ and in view of remark \ref{remarkmuisomega} when $p=2$ we have $\mu'_{\infty}=\omega'_{\infty}$ in \eqref{IpIn}. Moreover, it is shown in the proof of Theorem 4 in \cite{LNV11} that if $\theta_{0}=\theta_{0}(N,p)$ where $\theta_{0}$ is as in \eqref{nablaphi} is small enough (independent of $p$ when $p\in [3/2, 5/2]$) there exists $c=c(n)\geq 1$ and $\eta=\eta(\phi,n)$, $0<\eta<1/2$, such that
\begin{align}
\label{muomega}
c^{-1}\leq \frac{\tilde{I}(\phi,p)}{\tilde{I}(\phi,2)}\leq c\, \, \, \mbox{for}\, \, \, p\in(2-\eta, 2+\eta).
\end{align}
It follows from (\ref{IpIn}) and (\ref{muomega}) that for fixed $p\in(2-\eta, 2+\eta)$, p-harmonic measure $\mu'_{\infty}$ and harmonic measure $\omega'_{\infty}$ both have the Hausdorff dimension either $<n-1$ or $>n-1$ when $\theta_{0}$ is small enough. On the other hand,  Corollary 1 in \cite{LVV05} guarantees existence of a \textit{Wolff snowflake} $\tilde{\Omega}_{\infty}\subset\mathbb{R}^{n}$, $n\geq 3$, for which 
\begin{align}
\label{hdomegan-1}
I(\epsilon)<0 \, \, \implies \, \, \tilde{I}(\phi,2)<0
\end{align}
where $\omega'_{\infty}$ is harmonic measure for $\tilde{\Omega}_{\infty}$ with respect to a pole in the domain. Here we have used above observations and \eqref{IpIn} to obtain \eqref{hdomegan-1}. Let $p$ be fixed, $p\in(2-\eta, 2+\eta)$ and let $\tilde{u}_{\infty}$ be the p-harmonic function in $\tilde{\Omega}_{\infty}$ with zero continuous boundary values on $\partial\tilde{\Omega}_{\infty}$ and $|x_{n}-\tilde{u}_{\infty}(x)|\to 0$ uniformly as $|x|\to\infty$. Let $\tilde{\mu}_{\infty}$ be the p-harmonic measure associated with $\tilde{u}_{\infty}$ as in (\ref{pharmonic}). We first use \eqref{muomega} and then \eqref{IpIn} to get 
\[
\tilde{I}(\phi,2)<0\, \, \implies\, \, \hd{\tilde{\mu}_{\infty}}<n-1
\]
when $p\in(2-\eta, 2+\eta)$ for some $\eta>0$ small. Then we conclude that  Lemma \ref{dlessn-1} holds for $\tilde{\Omega}_{\infty}$ and $\tilde{\mu}_{\infty}$. We then repeat the argument in section \ref{proofofmain} to get the enlarged domain $\Omega_{\epsilon}^{+}$ from $\tilde{\Omega}_{\infty}$ as described in section \ref{proofofmain}. Let $\tilde{u}$ be a p-harmonic function in $\tilde{\Omega}:=\Omega_{\epsilon}^{+}$ with zero continuous boundary value and let $\tilde{\mu}$ be the p-harmonic measure associated with $\tilde{u}$ as in (\ref{pharmonic}). Following section \ref{proofofmain} we can find a compact set $\tilde{K}$ satisfying \eqref{muK}. in view of these observations we conclude the validity of Theorem \ref{main2}.
\qed
\begin{remark}
One of the main reason that Theorem \ref{main2} is stated for $n\geq 3$ is that our proof relies on Wolff's result in \cite{W95} which is valid only when $n\geq 3$.
\end{remark}
\section{Acknowledgments}
The author would like to thank Jonas Azzam for reading an earlier version of this manuscript and Matthew Badger for his suggestions. The author also thanks John Lewis for his suggestions and fruitful discussions to improve section \ref{proofofmain2}. The author has been supported in part by ICMAT Severo Ochoa project SEV-2011-0087. He acknowledges that the research leading to these results has received funding from the European Research Council under the European Union's Seventh Framework Programme (FP7/2007-2013)/ ERC agreement no. 615112 HAPDEGMT.

\def\cprime{$'$} \def\cprime{$'$}
\providecommand{\bysame}{\leavevmode\hbox to3em{\hrulefill}\thinspace}
\providecommand{\MR}{\relax\ifhmode\unskip\space\fi MR }
\providecommand{\MRhref}[2]{%
  \href{http://www.ams.org/mathscinet-getitem?mr=#1}{#2}
}
\providecommand{\href}[2]{#2}

\end{document}